\RequirePackage[l2tabu,orthodox]{nag}

\documentclass{scrartcl}

\KOMAoption{fontsize}{11pt}

\usepackage[T1]{fontenc} 
\usepackage[utf8]{inputenc}
\usepackage{lmodern} 
\usepackage{microtype}

\usepackage{amsmath}
\usepackage{amsthm}
\usepackage{amsfonts}
\usepackage{mathtools}

\usepackage[colorlinks=true,urlcolor=green,citecolor=green,linkcolor=green,linktocpage,pdfpagelabels,bookmarksnumbered,bookmarksopen]{hyperref}

\numberwithin{equation}{section}
\newtheorem{theorem}{Theorem}[section]
\newtheorem{lemma}[theorem]{Lemma}

\newtheorem{proposition}[theorem]{Proposition}

\theoremstyle{remark} 
\newtheorem{remark}[theorem]{Remark}

\theoremstyle{definition}

\begin{document}

\title{Existence of solutions for a semirelativistic Hartree equation with unbounded potentials}
\author{Simone Secchi \thanks{The author is supported by the MIUR 2015 PRIN project ``Variational methods, with applications to problems in mathematical physics and geometry''.}} 
\date{\today}
\maketitle

\hfil\emph{To Francesca, always}\hfill

\begin{abstract} 
	We prove the existence of a solution to the semirelativistic Hartree equation \[ \sqrt{-\Delta+m^2}u+ V(x) u = A(x)\left( W * |u|^p \right) |u|^{p-2}u \] under suitable growth assumption on the potential functions $V$ and $A$. In particular, both can be unbounded from above.
\end{abstract}


\bigskip

\section{Introduction}

The mean field limit of a quantum system describing many
self-gravitating, relativistic bosons with rest mass $m>0$ leads to
the time-dependent pseudo-relativistic Hartree equation
\begin{equation}\label{eq:1.1}
i\frac{\partial \psi}{\partial t} = \left(\sqrt{-\Delta +m^2} -m \right)\psi -  \left( \frac{1}{|x|} * |\psi|^2 \right)  \psi, \quad x \in \mathbb{R}^3
\end{equation}
where $\psi : \mathbb{R} \times \mathbb{R}^3 \to \mathbb{C}$ is the
wave field. Such a physical system is often referred to as a boson
star in astrophysics (see \cite{es,fl,fjl}).  Solitary wave solutions
$\psi(t,x)= e^{-i t \lambda} \phi$, $\lambda \in \mathbb{R}$ to
equation (\ref{eq:1.1}) satisfy the equation
\begin{equation*}
\left( \sqrt{-\Delta +m^2} -m \right)\phi - \left( \frac{1}{|x|} * |\phi|^2 \right)  \phi
=\lambda \phi \quad \mbox{in $\mathbb{R}^3$}.
\end{equation*}

For the non-relativistic Hartree equation driven by the \emph{local} differential operator $-\Delta+m^2$, existence and uniqueness (modulo translations) of a minimizer were proved by Lieb  \cite{lieb} by using symmetric decreasing rearrangement inequalities.
Within the same setting, always for the negative Laplacian, P.-L. Lions \cite{l2} proved existence of infinitely
many spherically symmetric solutions by application of abstract critical point theory both
with and without constraints for a more general radially symmetric convolution potential.
The  non-relativistic Hartree equations is also known as the Choquard-Pekard or Schr\"odinger-Newton equation and recently a large amount of papers are devoted to the study of solitary states and its semiclassical limit: see \cite{a,CSM,CT,cho,ccs,ccs1,css,CSProc,CNbis,CN,fl,ls,mpt,mz,pe2,pe3,s,t,ww} and references therein.


In this paper, which is somehow a continuation of the investigation we began in \cite{Secchi1,Secchi2}, we consider the equation
\begin{equation} \label{eq:1}
\sqrt{-\Delta+m^2}u+ V(x) u = A(x)\left( W * |u|^p \right) |u|^{p-2}u
\end{equation}
in the weighted space 
\[
X = \left\{ u \in L^2(\mathbb{R}^N) \mid \|u\|_X<\infty
\right\},
\]
where
\[
\|u\|_X^2 = \int_{\mathbb{R}^N} \left| \left(-\Delta+m^2 \right)^{1/4} u \right|^2 \, dx + \int_{\mathbb{R}^N} V(x) |u|^2 \, dx.
\]
In the case $A \equiv 1$, equation \eqref{eq:1} has been recently investigated  in \cite{CSProc}, where least-energy solutions are constructed for a bounded potential $V$.

\bigskip

We collect our assumptions.
\begin{itemize}
	\item[(H1)] $N >1$, $\beta>1$, $2 \leq p <  \frac{2N\beta}{\beta(N-1)}$ and $r > \max \left\{ 1,\frac{1-N\beta}{2+\beta N (p-2)-\beta p}\right\}$.
	\item[(H2)] $A \in L^\infty_{\mathrm{loc}}(\mathbb{R}^N)$ satisfies $A(x) \geq 1$ for almost every $x \in \mathbb{R}^N$.
	\item[(H3)] $V\colon \mathbb{R}^N \to \mathbb{R}$ is continuous and
	\begin{equation*} \label{ass:V}
	\inf_{\mathbb{R}^N} V > -B
	\end{equation*}
	for some $B>0$; furthermore,
	\[
	\lambda_1 = \inf_{u \in X \setminus \{0\}} \frac{\int_{\mathbb{R}^N} \left| (-\Delta+m^2)^{1/4}u \right|^2 + \int_{\mathbb{R}^N} V |u|^2}{\|u\|_2^2}>0.
	\]
	\item[(H4)] There exist $C_0>0$ and $R_0>0$ such that
	\[
	A(x)^\frac{2r}{2r-1} \leq C_0 \left( 1+ \left( \max\{0,V(x)\} \right)^\frac{1}{\beta} \right)
	\]
	for all $|x| > R_0$.
	\item[(H5)] $W \colon \mathbb{R}^N \setminus \{0\} \to \mathbb{R}$ and $W=W_1+W_2$, where $W_1 \in L^r$ and $W_2 \in L^\infty$.
	\end{itemize}
Equivalently, condition (H4) can be stated as
\begin{equation*}
\limsup_{|x| \rightarrow +\infty}\frac{A(x)^\frac{2r}{2r-1}}{1+ \left( \max \left\{ 0,V(x)\right\} \right)^\frac{1}{\beta}} < +\infty. 
\end{equation*}
\begin{remark}
	We highlight that both the
	potential $V$ and the potential $A$ in front of the nonlinear term in
	the right-hand side can be unbounded. To the best of our knowledge,
	this case is treated here for the first time. In the very recent paper \cite{CL}, a similar equation has been studied under an assumption that, in our framework, reads as $V \equiv 1$ and $\lim_{|x| \to +\infty} A(x)=1$.
\end{remark}

\medskip

To state our last assumption, we define, for any open subset $\Omega \subset \mathbb{R}^N$ and any $t \geq 2$, the quantity (see Section 2 for the definition of the space $L^{1/2,2}(\Omega)$)
\begin{align*}
\nu \left(t,\Omega \right) &= \inf_{\substack{u \in L^{1/2,2}(\Omega) \\ u \neq 0}} \frac{\int_\Omega \left| (-\Delta+m^2)^{1/4} u\right|^2 + \int_\Omega V |u|^2}{\left( \int_\Omega |u|^t \right)^{2/t}} \quad \text{if $\Omega \neq \emptyset$} \\
{}&= +\infty \quad\text{if $\Omega=\emptyset$}.
\end{align*}
We assume that
\begin{itemize}
\item[($\nu$)] there exists $t_0 \in [2,2N(N-1))$  such that
\[
\lim_{R \to +\infty} \nu(t_0,\mathbb{R}^N \setminus \overline{B_R}) = +\infty,
\]
where $B_R = \left\{ x \in \mathbb{R}^N \mid |x| <R \right\}$.
\end{itemize}
\bigskip 

The number~$\lambda_1$ in assumption (H3) can be actually seen as an \emph{eigenvalue},
as the following result shows.
\begin{proposition}
	Assume that $V \colon \mathbb{R}^N \to \mathbb{R}$ is continuous and
	$\inf_{\mathbb{R}^N} V > -B$ for some $B>0$. Assume also that
	($\nu$) holds with $t_0=2$. Then there exists $\varphi_1 \in X$ such
	that
	\begin{equation*}
	\sqrt{-\Delta+m^2}\varphi_1 + V(x) \varphi_1 = \lambda_1 \varphi_1 \quad \mbox{in $\mathbb{R}^N$.}
	\end{equation*}
	Furthermore $\varphi_1$ can be assumed to be positive in
	$\mathbb{R}^N$.
\end{proposition}
\begin{proof}
	It follows immediately from the assumptions that $\lambda_1 > -B
	>0$. Let $\{u_n\}_n$ be a minimizing sequence for $\lambda_1$, in
	the sense that
	\begin{align}
	\lim_{n \to +\infty} \int_{\mathbb{R}^N} \left| \left( -\Delta+m^2 \right)^\frac{1}{4} u_n \right|^2 + \int_{\mathbb{R}^N} V |u_n|^2 &= \lambda_1 \label{eq:1.5} \\
	\int_{\mathbb{R}^N} |u_n|^2 &=1. \label{eq:1.6}
	\end{align}
	As in \cite[Lemma 2.1]{ccs2}, we can assume that $u_n$ is
	non-negative. It follows from \eqref{eq:1.5} that
	\begin{equation*}
	\int_{\mathbb{R}^N} \left| (-\Delta+m^2)^\frac{1}{4} u_n \right|^2 \leq \lambda_1 + 1 + B.
	\end{equation*}
	Let $u$ be the weak limit of $\{u_n\}_n$ in
	$L^{1/2,2}(\mathbb{R}^N)$. Since $u_n \to u$ strongly in $L^2(B_R)$
	for any $R>0$, we may assume that $u_n$ converges to $u$ almost
	everywhere, and consequently $u \geq 0$. We fix a smooth cut-off
	function $\psi$ such that $\psi \equiv 0$ on $B_R$ and $\psi \equiv 1$
	on $\mathbb{R}^N \setminus B_{R+1}$. Then
	\begin{multline*}
	\left\| u_n - u \right\|_2^2 \leq \left\| (1-\psi) (u_n-u) \right\|_2^2 + \left\| \psi (u_n-u) \right\|_2^2 \\
	\leq o_n(1) + \frac{1}{\nu (2,\mathbb{R}^N \setminus B_R)} \left( \int_{\mathbb{R}^N} \left| \left( -\Delta+m^2 \right)^\frac{1}{4} (u_n-u) \right|^2 + \int_{\mathbb{R}^N} V |u_n-u|^2 \right) \\
	\leq o_n(1) + o_R(1),
	\end{multline*}
	where $\lim_{n \to +\infty} o_n(1)=0$ and $\lim_{R \to +\infty}
	o_R(1)=0$: this proves that $u_n \to u$ strongly in
	$L^2(\mathbb{R}^N)$. In particular, $\int_{\mathbb{R}^N} |u|^2
	=1$ due to \eqref{eq:1.6}. Plainly,
	\begin{equation*}
	\int_{\mathbb{R}^N} \left| (-\Delta+m^2)^\frac{1}{4} u \right|^2 \leq \liminf_{n \to +\infty} \left| (-\Delta+m^2)^\frac{1}{4} u_n \right|^2.
	\end{equation*}
	Set $G = \left\{ x \in \mathbb{R}^N \mid V(x)>1 \right\}$. Since $V
	\in L^\infty (\mathbb{R}^N \setminus G)$, we have
	\begin{equation*}
	\lim_{n \to +\infty} \int_{\mathbb{R}^N \setminus G} V |u_n|^2 = \int_{\mathbb{R}^N \setminus G} V |u|^2.
	\end{equation*}
	On the other hand, by \cite[Theorem 6.54]{FonsecaLeoni} we deduce that
	\begin{equation*}
	v \in L^2(G) \mapsto \int_{G} V |v|^2
	\end{equation*}
	is weakly lower semicontinuous. To summarize,
	\begin{align*}
	\liminf_{n \to +\infty} \int_{\mathbb{R}^N} V |u_n|^2 &= \lim_{n \to +\infty} \int_{\mathbb{R}^N \setminus G} V |u_n|^2 + \liminf_{n \to +\infty} \int_{\mathbb{R}^N} V |u_n|^2 \\
	&\geq \int_{\mathbb{R}^N} V |u|^2.
	\end{align*}
	Hence $u$ is a minimizer for $\lambda_1$. The strict positivity of $u$
	follows from \cite[Proposition 2]{FallFelli}.
\end{proof}
\begin{remark}
	In the previous Proposition, we did \emph{not} assume that
	$\lambda_1>0$. Under this additional assumption, the proof would be
	much easier.
\end{remark}

It is important to remark that condition ($\nu$) is related to some
other popular conditions. In Proposition \ref{prop:appendix1} we prove that
the coercivity assumption
\[
\lim_{|x| \to +\infty} V(x)=+\infty
\]
implies ($\nu$). Similarly, the condition
\begin{itemize}
\item[($\nu'$)] For every $b>0$, the set $V^b$ has finite Lebesgue
  measure
\end{itemize}
also implies ($\nu$), see Proposition \ref{prop:appendix2}. Finally,
Sirakov's condition \cite{Sirakov}
\begin{itemize}
	\item[($\nu''$)] There exists $t_1 \in [2,2N/(N-1))$ such that, for any $r>0$ and any sequence $\{x_n\}_n$ of points in $\mathbb{R}^N$ such that $\lim_{n \to +\infty} |x_n| = +\infty$, there results $\lim_{n \to +\infty} \nu(t_1,B(x_n,r))=+\infty$
\end{itemize}
is equivalent to ($\nu$): see Proposition
\ref{prop:appendix3}. 

\bigskip

\begin{remark}
	While dealing with the semirelativistic Hartree equation, it is customary to rewrite the operator $\sqrt{-\Delta+m^2}+V$ as $\sqrt{-\Delta+m^2}-m+(V+m)$ in order to exploit the fact that $\sqrt{-\Delta+m^2}-m>0$ in the sense of functional calculus. The natural assumption from a variational point of view is therefore $\inf_{\mathbb{R}^N} V > -m$. Our assumption (H3) is, in general, less restrictive as it only requires  a \emph{suitable} lower bound for $V$.
\end{remark}

\medskip
We can state the main result of this paper.
\begin{theorem} \label{th:main}
	Suppose that (H1), (H2), (H3), (H4), (H5), and ($\nu$) are satisfied. Then equation \eqref{eq:1} has infinitely many distinct solutions.
\end{theorem}
The proof is based on the ideas developed in \cite{Sirakov} for a local equation and a pointwise nonlinearity.
\begin{remark}
	It should be noticed that our results continue to hold if we replace \eqref{eq:1} with
	\begin{equation*}
	(-\Delta+m^2)^{s} u + V(x) u = \left(W*|u|^{p} \right)|u|^{p-2}u
	\end{equation*}
	with $0<s<1$ and $N > 2s$, for instance $N \geq 3$. Of course some numbers must be replaced: $2N/(N-1)$ should become $2N/(N-2s)$, and so on. We prefer to work out the details for $s=1/2$, which corresponds to the physical model of the Hartree equation.
\end{remark}
In Section 2 we introduce the necessary preliminaries on function space. In Section 3 we compare several assumptions the ensure a compact embedding result. In Section 4 we prove a compactness theorem that is used in Section 5 to prove our main existence result.

\bigskip

\begin{minipage}{4in} \textbf{Notation}
	\begin{enumerate}
		\item The letters $c$ and $C$ will stand for a generic positive
		constant that may vary from line to line.
		\item The operator $d$ will be reserved for the (Fr\'{e}chet)
		derivative, so that $dI$ denotes the Fr\'{e}chet derivative of a function $I$.
		\item The symbol $\mathcal{L}^N$ will be reserved for the
		Lebesgue $N$-dimensional measure.
		\item The Fourier transform of a function $u$ will be denoted
		by $\mathcal{F}u$.
	\end{enumerate}
\end{minipage}


%
\section{Variational setting}
Let us recall the definition of the \emph{Bessel function space} 
defined for~$\alpha>0$ by
\[ 
L^{\alpha,2}(\mathbb{R}^N) = \left\{ f \colon f=G_\alpha \star g \
\text{for some $g \in L^2(\mathbb{R}^N)$} \right\},
\] 
where the Bessel convolution kernel is defined by
\begin{equation*} 
G_\alpha (x) = \frac{1}{(4 \pi )^{\alpha
		/2}\Gamma(\alpha/2)} \int_0^\infty \exp \left( -\frac{\pi}{t} |x|^2
\right) \exp \left( -\frac{t}{4\pi} \right) t^{\frac{\alpha - N}{2}-1}
\, dt
\end{equation*} 
The norm of this Bessel space is $\|f\| = \|g\|_2$ if
$f=G_\alpha \star g$. The operator 
\[
(I-\Delta)^{-\alpha} u = G_{2\alpha}
\star u
\] 
is usually called Bessel operator of order $\alpha$.
In Fourier variables the same operator reads
\begin{equation*} 
G_\alpha = \mathcal{F}^{-1} \circ
\left( \left(1+|\xi|^2 \right)^{-\alpha /2} \circ \mathcal{F} \right),
\end{equation*} 
so that
\begin{equation*} 
\|f\| = \left\| (I-\Delta)^{\alpha /2} f \right\|_2.
\end{equation*}
For more detailed information, see \cite{Adams, Stein} and the
references therein. The use of $(-\Delta+m^2)^{\alpha /2}$ instead of $(-\Delta+I)^{\alpha/2}$ is clearly harmless. We summarize the embedding properties of Bessel spaces. For the proofs
we refer to \cite[Theorem 3.2]{Felmer}, \cite[Chapter V, Section 3]{Stein} and \cite[Section 4]{Strichartz}.

\begin{theorem} \label{th:bessel}
	\begin{enumerate}
		\item $L^{\alpha,2}(\mathbb{R}^N) =
		W^{\alpha,2}(\mathbb{R}^N) = H^\alpha (\mathbb{R}^N)$.
		\item If $\alpha \geq 0$ and $2 \leq q \leq
		2_\alpha^*=2N/(N-2\alpha)$, then $L^{\alpha,2}(\mathbb{R}^N)$ is
		continuously embedded into $L^q(\mathbb{R}^N)$; if $2 \leq q <
		2_\alpha^*$ then the embedding is locally compact.
		\item Assume that $0 \leq \alpha \leq 2$ and $\alpha >
		N/2$. If $\alpha -N/2 >1$ and $0< \mu \leq \alpha - N/2-1$, then
		$L^{\alpha,2}(\mathbb{R}^N)$ is continuously embedded into
		$C^{1,\mu}(\mathbb{R}^N)$. If $\alpha -N/2 <1$ and $0 < \mu \leq
		\alpha -N/2$, then $L^{\alpha,2}(\mathbb{R}^N)$ is continuously
		embedded into $C^{0,\mu}(\mathbb{R}^N)$.
	\end{enumerate}
\end{theorem}
\begin{remark} 
	Although the Bessel space $L^{\alpha,2}(\mathbb{R}^N)$
	is topologically undistinguishable from the So\-bo\-lev fractional space
	$H^\alpha(\mathbb{R}^N)$, we will not systematically confuse them, since our equation
	involves the Bessel norm.
\end{remark}
For a general open subset $\Omega$ of $\mathbb{R}^N$, the Bessel space $L^{\alpha,2}(\Omega)$
is defined ad the space of the restrictions to $\Omega$ of functions in $L^{\alpha,2}(\mathbb{R}^N)$. In this paper we will always take $\alpha = 1/2$.
As we said in the Introduction, we work in the weighted space 
\[
X = \left\{ u \in L^{1/2,2}(\mathbb{R}^N) \mid \|u\|_X<\infty
\right\}
\]
where
\[
\|u\|_X^2 = \int_{\mathbb{R}^N} \left| \left(-\Delta+m^2 \right)^{1/4} u \right|^2 \, dx + \int_{\mathbb{R}^N} V(x) |u|^2 \, dx
\]
\begin{lemma} \label{lem:2.1}
If assumption (H3) holds, then there exists a constant   $\aleph>0$ such that
\[
\|u\|_X\geq \aleph \|u\|_{L^{1/2,2}(\mathbb{R}^N)}
\]
for every $u \in X$.
\end{lemma}
\begin{proof}
We proceed by contradiction. Assume that for some sequence $\{u_n\}_n$ in $X$ there results
\[
\|u_n\|_{L^{1/2,2}}^2=1, \quad \|u_n\|_X^2 < \frac{1}{n}.
\]
By definition of $\lambda_1>0$, the last inequality entails $\|u_n\|_{L^2} =o(1)$. But then the contradiction
\[
o(1) =-B \|u_n\|_{L^2}^2 \leq \int_{\mathbb{R}^N} V |u_n|^2 < \frac{1}{n}-1
\]
arises as $n \rightarrow +\infty$.
\end{proof}
Solutions to \eqref{eq:1} correspond to critical points of the functional $I \colon X \rightarrow \mathbb{R}$ defined by
\begin{equation} \label{eq:I}
I(u) = \frac{1}{2} \|u\|_X^2 - \frac{1}{2p} \int_{\mathbb{R}^N \times \mathbb{R}^N} W(x-y)A(x)|u(x)|^p |u(y)|^p \, dx \, dy.
\end{equation}
We need to prove that $I$ is well-defined on $X$.
\begin{proposition} \label{prop:embedding}
	The space $X$ is continuously embedded into the weighted Lebesgue space \[
	L^t(\mathbb{R}^N,A^\frac{2r}{2r-1}d\mathcal{L}^N)
	\]
	for every $2 \leq t \leq \frac{2(N\beta-1)}{\beta(N-1)}$.
\end{proposition}
\begin{proof}
Let us decompose
\begin{multline}\label{eq:2.2.1}
\int_{\mathbb{R}^N} A(x)^\frac{2r}{2r-1}|u(x)|^t \, dx = \int_{B(0,R_0)} A(x)^\frac{2r}{2r-1}|u(x)|^t \, dx \\
{}+ \int_{\mathbb{R}^N \setminus B(0,R_0)} A(x)^\frac{2r}{2r-1}|u(x)|^t \, dx,
\end{multline}
where $R_0>0$ is the number defined in assumption (H4). Now,
\begin{multline*}
\int_{\mathbb{R}^N \setminus B(0,R_0)} A(x)^\frac{2r}{2r-1}|u(x)|^t \, dx \leq C_0 \int_{\mathbb{R}^N \setminus B(0,R_0)} \left( 1+ \left( \max \{0,V(x)\} \right)^\frac{1}{\beta} \right) |u(x)|^t \, dx \\
\leq C_0 \int_{\Omega_1} |u(x)|^t \, dx + C_0 \int_{\Omega_2} \left( 1+ V(x)^\frac{1}{\beta} \right)|u(x)|^t \, dx,
\end{multline*}
where
\begin{align*}
\Omega_1 = \left\{ x \in \mathbb{R}^N \mid V(x)<0\right\}, \qquad
\Omega_2 = \left\{ x \in \mathbb{R}^N \mid V(x)\geq 0\right\}.
\end{align*}
As a consequence, inserting this into \eqref{eq:2.2.1},
\[
\int_{\mathbb{R}^N} A(x)^\frac{2r}{2r-1}|u(x)|^t \, dx \leq 2C_0 \|u\|_{L^t}^t + \|A\|_{L^\infty(B(0,R_0))} \|u\|_{L^t}^t + \int_{\Omega_2} V(x)^\frac{1}{\beta} |u(x)|^t \, dx.
\]
If we apply H\"{o}lder's inequality to the last integral, we obtain
\begin{multline*}
\int_{\Omega_2} V(x)^\frac{1}{\beta} |u(x)|^t \, dx = \int_{\Omega_2} V(x)^\frac{1}{\beta} |u(x)|^\frac{2}{\beta} |u(x)|^{t-\frac{2}{\beta}}\, dx \\
\leq \left( \int_{\Omega_2} V(x) |u(x)|^2 \, dx \right)^\frac{1}{\beta} \left( \int_{\Omega_2} |u(x)|^{\left(t-\frac{2}{\beta} \right) \frac{\beta}{\beta-1}} \, dx \right)^\frac{\beta-1}{\beta}.
\end{multline*}
But
\begin{multline*}
\int_{\Omega_2} V(x) |u(x)|^2 \, dx = \int_{\mathbb{R}^N} V(x) |u(x)|^2 \, dx - \int_{\left\{x \in \mathbb{R}^N \mid V(x) < 0\right\}} V(x) |u(x)|^2 \, dx \\
\leq \|u\|_X^2 + B \int_{\mathbb{R}^N} |u(x)|^2 \, dx
\leq \left(1+\frac{B}{\lambda_1} \right) \|u\|_X^2.
\end{multline*}
In conclusion,
\begin{equation} \label{eq:2.2}
\int_{\mathbb{R}^N} A(x)^\frac{2r}{2r-1}|u(x)|^t \, dx \leq C_1 \|u\|_{L^t}^t + C_2 \|u\|_X^\frac{2}{\beta} \|u\|_{L^{(t\beta -2)/(\beta-1)})}^{(t\beta-2)/\beta}.
\end{equation}
It is elementary to check that
\begin{equation*} \label{eq:2.3}
2 \leq \frac{t\beta -2}{\beta-1} \leq \frac{2N}{N-1}
\end{equation*}
whenever 
\begin{equation*} \label{eq:2.4}
2 \leq t \leq \frac{2(N\beta-1)}{\beta(N-1)},
\end{equation*}
so that we can invoke the continuous embedding of $X$ into $L^s(\mathbb{R}^N)$ for $2 \leq s \leq 2N/(N-1)$ and conclude from \eqref{eq:2.2} that there exists a positive constant $C$ such that
\[
\int_{\mathbb{R}^N} A(x)^\frac{2r}{2r-1} |u(x)|^t \, dx \leq C \|u\|_X^t \quad\text{for every $u \in X$}.
\]
This completes the proof.
\end{proof}	
\begin{remark}
	Observe that when $\inf_{\mathbb{R}^N} V>0$ and $A =1$ identically, we can let $\beta \to +\infty$ and recover the weaker assumption $2 \leq t \leq 2N/(N-1)$. 
\end{remark}
To proceed further, we need the following inequality due to Hardy, Littlewood and Sobolev. We firstly recall that a function $h$ belongs to the weak $L^q$ space $L_w^q(\mathbb{R}^N)$ if there exists a constant $C>0$ such that, for all $t>0$, 
\[
\mathcal{L}^N \left( \left\{ x \in \mathbb{R}^N \mid |h(x)| > t \right\} \right) \leq \frac{C^q}{t^q}.
\]
Its norm is then
\[
\|h\|_{q,w} = \sup_{t>0} t \left(\mathcal{L}^N \left( \left\{ x \in \mathbb{R}^N \mid |h(x)| > t \right\} \right) \right)^{1/q}.
\]
\begin{proposition}[\cite{Lieb1983}]
	Assume that $p$, $q$ and $t$ lie in $(1,+\infty)$ and $p^{-1}+q^{-1}+t^{-1}=2$. Then, for some constant $N_{p,q,t}>0$ and for any $f \in L^p(\mathbb{R}^N)$, $g \in L^t(\mathbb{R}^N)$ and $h \in L_w^q(\mathbb{R}^N)$, we have the inequality
	\begin{equation} \label{eq:hls}
	\left| \int_{\mathbb{R}^N \times \mathbb{R}^N} f(x)h(x-y)g(y)\, dx\, dy
	\right| \leq N_{p,q,t} \|f\|_p \|g\|_t \|h\|_{q,w}.
	\end{equation}
\end{proposition}
Writing $W=W_1+W_2 \in L^r(\mathbb{R}^N) + L^\infty(\mathbb{R}^N)$  and using $(\ref{eq:hls})$ we
can estimate the convolution term as follows by means of Proposition \ref{prop:embedding}:
\begin{align} \label{eq:15}
\int_{\mathbb{R}^N} (W *|u|^p)A|u|^p &= \int_{\mathbb{R}^N\times \mathbb{R}^N} |u(x)|^p A(x) W(x-y) |u(y)|^p \, dx \, dy \nonumber \\
&= \int_{\mathbb{R}^N\times \mathbb{R}^N} |u(x)|^p A(x) W_1(x-y) |u(y)|^p \, dx \, dy \nonumber \\
&\qquad {} + \int_{\mathbb{R}^N\times \mathbb{R}^N} |u(x)|^p A(x) W_2(x-y)|u(y)|^p \, dx \, dy \nonumber \\
&\leq C \left( \int_{\mathbb{R}^N} |u|^{\frac{2rp}{2r-1}} A^{\frac{2r}{2r-1}} \right)^\frac{2r-1}{2r} \|W_1\|_r \left( \int_{\mathbb{R}^N} |u|^{\frac{2rp}{2r-1}}  \right)^\frac{2r-1}{2r} \nonumber \\
&\qquad {} + \|W_2\|_\infty \left( \int_{\mathbb{R}^N} |u(x)|^{p} A(x) \, dx \right) \left( \int_{\mathbb{R}^N} |u(y)|^{p}\, dy  \right) \nonumber \\
&\leq C \left( \int_{\mathbb{R}^N} |u|^{\frac{2rp}{2r-1}} A^{\frac{2r}{2r-1}} \right)^\frac{2r-1}{2r} \|W_1\|_r \left( \int_{\mathbb{R}^N} A^{\frac{2r}{2r-1}}|u|^{\frac{2rp}{2r-1}}  \right)^\frac{2r-1}{2r} \nonumber \\
&\qquad{}+ \|W_2\|_\infty \left( \int_{\mathbb{R}^N} |u|^{p} A^\frac{2r}{2r-1} \right) \left( \int_{\mathbb{R}^N} |u|^{p}A^\frac{2r}{2r-1}  \right) \nonumber \\
&\leq C \|u\|_{L^\frac{2rp}{2r-1}(\mathbb{R}^N,A^\frac{2r}{2r-1}d\mathcal{L}^N)}^{2p} \|W_1\|_r + \|W_2\|_\infty \|u\|_{L^p(\mathbb{R}^N,A^\frac{2r}{2r-1}d\mathcal{L}^N)}^{2p},
\end{align}
where we have used several times the fact that $A(x)\geq 1$ for almost every $x \in \mathbb{R}^N$.
Since
\[
r > \max \left\{ 1,\frac{1-N\beta}{2+\beta N(p-2)-\beta p} \right\} \quad\text{and}\quad p < \frac{2N\beta}{\beta(N-1)},
\]
we can use Proposition \ref{prop:embedding} and from (\ref{eq:15}) we see that the convolution term in $I$ is finite. It is easy to check, by the same token and taking into account the assumption $p\geq 2$, that
$I \in C^1(X)$.
\begin{remark}
	As before, if $\inf_{\mathbb{R}^N} V>0$ and $A=1$ identically, we can recover the assumption \[ r> \max \left\{ 1,\frac{N}{p+(2-p)N}\right\} \] used in \cite{CSProc}.
\end{remark}

\section{Comparison between conditions on \(V\)}

In this section we prove that condition ($\nu$) is actually weaker than both the coercivity of $V$ and of condition ($\nu'$).
We start with a preliminary technical result.
\begin{lemma} \label{lem:3.1}
	Let $\Omega \subset \mathbb{R}^N$ be an open set. If $2 \leq t < r< 2N/(N-1)$, then there exist a constant $C=C(r,t,N,\lambda_1)>0$ and a number $0<\theta_1<1$ such that
	\begin{equation*}
	\nu(r,\Omega) \geq C \left( \nu(t,\Omega) \right)^{\theta_1}.
	\end{equation*}
	If $2<r<t<2N/(N-1)$, then there exist a constant $C=C(r,t,N,\lambda_1)>0$ and a number $0<\theta_2<1$ such that
	\begin{equation*}
	\nu(r,\Omega) \geq C \left( \nu(t,\Omega) \right)^{\theta_2}.
	\end{equation*}
	In particular,
	\[
	\lim_{R \to +\infty} \nu(r,\mathbb{R}^N \setminus \overline{B_R}) = +\infty
	\]
	if
	\[
	\lim_{R \to +\infty} \nu(t,\mathbb{R}^N \setminus \overline{B_R}) = +\infty.
	\]
\end{lemma}
\begin{proof}
	We write
	\[
	\frac{1}{r} = \frac{1-\theta}{t} + \frac{\theta}{\frac{2N}{N-1}}
	\]
	for some $0<\theta<1$. Then $\|u\|_r \leq \|u\|_t^{1-\theta} \|u\|_{\frac{2N}{N-1}}^\theta$ for every $u \in L^{1/2,2}(\Omega)$. The Gagliardo-Nirenberg inequality (see \cite[Theorem 2.1]{Park} together with Theorem \ref{th:bessel}) yields $\|u\|_r \leq C \|u\|_{L^{1/2,2}}^\theta \|u\|_t^{1-\theta}$ for every $u \in L^{1/2,2}(\Omega)$. We invoke Lemma \ref{lem:2.1} to get
	\[
	\|u\|_r^2 \leq C \|u\|_{X}^{2\theta} \|u\|_t^{2(1-\theta)}
	\]
	As a consequence,
	\begin{multline*}
	\nu(r,\Omega) = \inf_{u \in L^{1/2,2}(\Omega) \setminus \{0\}} \frac{\|u\|_X^2}{\|u\|_r^2} \geq \frac{1}{C} \inf_{u \in L^{1/2,2}(\Omega) \setminus \{0\}}  \frac{\|u\|_X^2}{\|u\|_{V}^{2\theta} \|u\|_t^{2(1-\theta)}} \\
	\geq \frac{1}{C} \inf_{u \in L^{1/2,2}(\Omega) \setminus \{0\}} \frac{\|u\|_X^{2(1-\theta)}}{\|u\|_t^{2(1-\theta)}} 
	\geq \frac{1}{C} \left( \nu(t,\Omega) \right)^{1-\theta}.
	\end{multline*}
	The second inequality follows by the same token.
\end{proof}
\begin{proposition} \label{prop:appendix1}
	If $\lim_{|x| \rightarrow +\infty} V(x)=+\infty$, then ($\nu$) holds true.
\end{proposition}
\begin{proof}
	By Lemma \ref{lem:3.1} it is sufficient to prove the validity of $(\nu$) with $t_0=2$. For any $u \in X$ and any $R>0$, 
	\begin{align*}
		\int_{\mathbb{R}^N \setminus B_R} \left| (-\Delta+m^2)^\frac{1}{4} u \right|^2 + \int_{\mathbb{R}^N \setminus B_R} V |u|^2 \geq \int_{\mathbb{R}^N \setminus B_R} V |u|^2 \geq \inf_{\mathbb{R}^N\setminus B_R} V \int_{\mathbb{R}^N \setminus B_R} |u|^2. 
	\end{align*}
	Therefore $\nu(2,\mathbb{R}^N \setminus B_R) \geq \inf_{\mathbb{R}^N\setminus B_R} V$, and we conclude by letting $R \to +\infty$.
\end{proof}
\begin{lemma} \label{lem:5.2}
	Let $\{\omega_n\}_n$ be a sequence of open subset of $\mathbb{R}^N$ such that
	\[
	\lim_{n \to +\infty} \mathcal{L}^N(\omega_n) = 0.
	\]
	Then, for every $\varrho>0$ and every $2 \leq t < 2N/(N-1)$, there results
	\begin{equation*}
	\lim_{n \to +\infty} \sup \left\{  \int_{\omega_n} |u|^t \mid \|u\|_{L^{1/2,2}} \leq \varrho \right\} =0. 
	\end{equation*}
\end{lemma}
\begin{proof}
	Since
	\[
	\int_{\omega_n} |u|^t \leq \|u\|_{L^{2N/(N-1)}}^t \mathcal{L}^N (\omega_n)^\frac{(2-t)N+t}{2N},
	\]
	by the Sobolev embedding we can find a constant $C_1>0$ such that
	\[
	\int_{\omega_n} |u|^t \leq C_1 \|u\|_{L^{1/2,2}}^t \mathcal{L}^N (\omega_n)^\frac{(2-t)N+t}{2N}
	\leq C_2 \varrho^t \mathcal{L}^N (\omega_n)^\frac{(2-t)N+t}{2N}
	\]
	whenever $\|u\|_{L^{1/2,2}} \leq \varrho$. The conclusion follows immediately.
\end{proof}
\begin{proposition} \label{prop:appendix2}
	Condition ($\nu'$) implies condition ($\nu$).
\end{proposition}
\begin{proof}
	For the sake of contradiction, let us suppose that there exist a sequence $R_n \to +\infty$ of positive numbers and a sequence $\{u_n\}_n$ of functions from $X$ such that
	\begin{align}
		&\sup_{n\in\mathbb{N}} \int_{\mathbb{R}^N \setminus B_{R_n}} \left| (-\Delta+m^2)^\frac{1}{4}u_n \right|^2 + \int_{\mathbb{R}^N \setminus B_{R_n}} V |u_n|^2 < +\infty \nonumber \\
		&\int_{\mathbb{R}^N\setminus B_{R_n}} |u_n|^2 =1. \label{eq:5.1}
	\end{align}
	By Lemma \ref{lem:2.1}, the sequence $\{u_n\}$ is bounded in $L^{1/2,2}(\mathbb{R}^N)$. For any $M>0$ we introduce the set
	\[
	\Omega(M,n) = \left\{ x \in \mathbb{R}^N \setminus B_{R_n} \mid V(x) < M\right\},
	\]
	so that $\lim_{n \to +\infty} \mathcal{L}^N(\Omega(M,n)) =0$. For some $C>0$, 
	\begin{multline} \label{eq:5.3}
		C \geq \int_{B_{R_n} \setminus \Omega(M,n)} V|u_n|^2 + \int_{\Omega(M,n)} V |u_n|^2 \geq \\
		M \int_{B_{R_n}} |u_n|^2 -(B+M) \int_{\Omega(M,n)} |u_n|^2.
	\end{multline}
	By Lemma \ref{lem:5.2} the last term of \eqref{eq:5.3} converges to zero as $n \to +\infty$. Since $M>0$ is arbitrary and \eqref{eq:5.1} holds true, we derive a contradiction.
\end{proof}
Finally, we prove that our assumption ($\nu$) is logically equivalent to Sirakov's condition ($\nu''$).
\begin{proposition} \label{prop:appendix3}
	Condition ($\nu''$) is equivalent to ($\nu$).
\end{proposition}
\begin{proof}
	\textbf{Step 1.} Assume that ($\nu$) holds. Given any $r>0$, we know that
	\[
	\nu(t_0,\mathbb{R}^N \setminus B_{R_n}) \leq \nu (t_0,B(x_n,r))
	\]
	where $R_n = \frac{1}{2}|x_n|-r$. Hence also ($\nu''$) holds.
	
	\noindent\textbf{Step 2.} On the contrary, assume that ($\nu$) does \emph{not} hold. Hence there are sequences $\{u_n\}_n$ in $X$ and $\{R_n\}_n$ in $(0,+\infty)$ such that $u_n \in L^{1/2,2}(\mathbb{R}^N)$, $\operatorname{supp}u_n \subset \mathbb{R}^N \setminus \overline{B(0,R_n)}$, $\lim_{n \to +\infty} R_n = +\infty$,
	\begin{align*}
		\int_{\mathbb{R}^N} \left| (-\Delta+m^2)^\frac{1}{4}u_n \right|^2 + \int_{\mathbb{R}^N} V|u_n|^2 &\leq C \\
		\int_{\mathbb{R}^N} |u_n|^{t_0} &=1.
	\end{align*}
	As before, the sequence $\{u_n\}_n$ is bounded in $L^{1/2,2}(\mathbb{R}^N)$ and by a Lemma of P.-L. Lions (see \cite[Lemma 2.4]{Secchi-ground}) there exist a sequence $\{x_n\}_n$ of points in $\mathbb{R}^N$, two constants $r_0>0$, $C_0>0$, such that, up to a subsequence,
	\begin{equation} \label{eq:5.3bis}
		\int_{B(x_n,r_0)} |u|^{t_0} \geq C_0.
	\end{equation}
	It follows from these properties that $|x_n| \to +\infty$. We fix a sequence of smooth cut-off functions $\varphi_n$ such that $\varphi_n \in C_c^\infty(\mathbb{R}^N)$, $0 \leq \varphi_n \leq 1$
	\begin{align*}
		\varphi(x) &= 1 \quad \text{if $x \in B(x_n,r_0)$}\\
		&= 0 \quad \text{if $x \notin B(x_n,2r_0)$}.
	\end{align*}	
	Let $v_n = \varphi_n u_n \in L^{1/2,2}(B(x_n,2r_0))$, so that by \eqref{eq:5.3bis}
	\[
	\int_{B(x_n,2r_0)} |v_n|^{t_0} \geq C_0.
	\]
	It is well known, see \cite{DPV}, that $\|v_n\|_{H^{1/2}} \leq C_1 \|u_n\|_{H^{1/2}}$ for some constant $C_1>0$. By Theorem \ref{th:bessel} and Lemma \ref{lem:2.1}, $\|v_n\|_X \leq C_2 \|u_n\|_X$. Thus for any $n \in \mathbb{N}$,
	\begin{equation*}
		\nu(t_0,B(x_n,2r_0)) \leq \frac{\int_{B(x_n,2r_0)} \left| (-\Delta+m^2)^\frac{1}{4} u_n \right|^2 + \int_{B(x_n,2r_0)}V|u_n|^2}{\left( \int_{B(x_n,2r_0)} |u_n|^{t_0} \right)^{2/t_0}}
		\leq \frac{C_2 \|u_n\|_X^2}{C_0^{2/t_0}} \leq C_3.
	\end{equation*}
	We have proved that also ($\nu''$) does \emph{not} hold.
\end{proof}
\section{A compact embedding theorem}
Theorem \ref{th:main} will be proved by applying the Symmetric Mountain Pass Theorem of Ambrosetti and Rabinowitz \cite{ar} to the functional $I$. The required compactness is recovered by embedding the space $X$ into a weighted Lebesgue space, see Proposition \ref{prop:3.2}.

We are ready to prove the main compactness result of this paper.
\begin{proposition} \label{prop:3.2}
Assume that (H1), (H2), (H3), (H4) and ($\nu$) are satisfied. Then $X$ is compactly embedded into $L^t(\mathbb{R}^N)$ for all $2 \leq t \leq \frac{2(N\beta-1)}{\beta(N-1)}$, and compactly embedded into $L^t(\mathbb{R}^N,A^{2r/(2r-1)}\, d\mathcal{L}^N)$ if $2 \leq t < \frac{2(N\beta-1)}{\beta(N-1)}$.
\end{proposition}
\begin{proof}
Consider any bounded sequence $\{u_n\}_n$ in $X$. By reflexivity, we assume without loss of generality that $u_n \rightharpoonup u$ weakly in $X$ as $n \to +\infty$. Choose a smooth function $\varphi \colon \mathbb{R}^N \to [0,1]$ with the property that $\varphi(x) = 0$ if $|x| \leq R$, while $\varphi(x) =1$ if $|x|>R+1$. Then
\[
\|u_n-u\|_t \leq \|(1-\varphi) (u_n-u) \|_{L^t(B_{R+1})} + \| \varphi (u_n-u)\|_{L^t(\mathbb{R}^N \setminus B_{R})}.
\]
The space $L^{1/2,2}(B_{R+1})$ is compactly embedded into $L^t(B_{R+1})$ because $\frac{2(N\beta-1)}{\beta(N-1)} < 2N/(N-1)$, see Theorem \ref{th:bessel}. Hence, up to a subsequence,
\[
\lim_{n \to +\infty} \|(1-\varphi) (u_n-u) \|_{L^t(B_{R+1})} =0.
\]
But by definition of $\nu$
\[
\nu(t,\mathbb{R}^N \setminus \overline{B_R}) \leq \frac{\|\varphi (u_n-u)\|_X^2}{\| \varphi (u_n-u)\|_{L^t(\mathbb{R}^N \setminus B_R)}}.
\]
We deduce that, using the fact that $\|\varphi (u_n-u)\|_{H^{1/2}} \leq C \|u_n-u\|_{H^{1/2}}$ and that $L^{1/2,2}=H^{1/2}$,
\[
\|\varphi (u_n-u) \|_{L^t(\mathbb{R}^N \setminus B_R)} \leq \frac{C}{\nu(t,\mathbb{R}^N \setminus \overline{B_R})},
\]
and we conclude by Lemma \ref{lem:3.1}. The second part follows directly from Proposition \ref{prop:embedding} (see in particular \eqref{eq:2.2}).
\end{proof}
\section{Existence of critical points}
In this section we apply the celebrated Mountain Pass Theorem and the embedding result proved earlier to find a critical point of the functional $I$ defined in \eqref{eq:I}.
\begin{proposition} \label{prop:geometry}
	The functional $I$ has the Mountain Pass geometry.
\end{proposition}
\begin{proof}
	\noindent\textbf{Step 1:} there exists a number $\rho>0$ such that $\|u\|_X = \rho$ implies $I(u)>0$.
	
	Indeed, by \eqref{eq:15} we find that for some constant $C>0$
	\[
	I(u) \geq \frac{1}{2}\|u\|_X^2 - C \|u\|_X^{2p}.
	\]
	Hence $I(u)\leq 0$ implies $\|u\|_X \geq C$ for another constant $C>0$. The claim follows easily.
	
	\noindent\textbf{Step 2:} there exists $e \in X$ such that $\|e\|_X > \rho$ and $I(e)<0$.
	
	Simply, pick $u \in X$ such that $\int_{\mathbb{R}^N} (K*|u|^{2p})|u|^{2p}>0$ and compute
	\[
	\lim_{t \to +\infty} I(tu) = \lim_{t \to +\infty} \frac{t^2}{2} \|u\|_X^2 - \frac{t^{2p}}{2p} \int_{\mathbb{R}^N} (W*|u|^{2p})|u|^{2p} = -\infty.
	\]
\end{proof}
In order to apply the Mountain Pass Theorem, we need to ensure the validity of the Palais-Smale compactness condition. This follows from Proposition \ref{prop:3.2}, as we now show.
\begin{proposition} \label{prop:PS}
	If $\{u_n\}_n$ is any sequence from $X$ such that $dI(u_n) =o(1)$ in $X^*$ and $I(u_n) \leq C$ for every $n \in \mathbb{N}$, then $\{u_n\}_n$ converges strongly in $X$ up to a subsequence.
\end{proposition}
\begin{proof}
First of all, the sequence 	$\{u_n\}_n$ is bounded in $X$. Indeed,
\[
C+o(1) \geq I(u_n)-\frac{1}{2p} dI(u_n)(u_n) = \left( \frac{1}{2} - \frac{1}{2p} \right) \|u_n\|_X^2.
\]
By reflexivity, we assume without loss of generality that $u_n \rightharpoonup u$ weakly in $X$ as $n \to +\infty$. It is plain that $u$ is a critical point of $I$. Now Proposition \ref{prop:3.2} implies that --- up to a subsequence --- $u_n \to u$ strongly in $L^{\frac{2rp}{2r-1}}(\mathbb{R}^N,|A|^{\frac{2r}{2r-1}}\, d\mathcal{L}^N)$, so that \eqref{eq:15} easily implies
\begin{equation*}
\lim_{n \to +\infty} \|u_n\|_X^2 = \lim_{n \rightarrow +\infty} \int_{\mathbb{R}^N}  (W*|u_n|^{2p})A|u_n|^{2p} = \int_{\mathbb{R}^N}  (W*|u|^{2p})A|u|^{2p} = \|u\|_X^2.
\end{equation*}
We conclude that $u_n \to u$ strongly in $X$ as $n \to +\infty$, and the proof is complete.
\end{proof} 
\begin{proof}[Proof of Theorem \ref{th:main}]
Let us remark that the right-hand side of equation \eqref{eq:1} is odd with respect to $u$. By Proposition \ref{prop:geometry} and Proposition \ref{prop:PS}, we can invoke the Symmetric Mountain Pass Theorem \cite{ar} and conclude that the functional $I$ possesses infinitely many critical points in $X$.
\end{proof}

\bigskip

\begin{minipage}{4in}
Simone Secchi \\
Dipartimento di Matematica e
Applicazioni \\ 
Universit\`a degli Studi di Milano Bicocca \\ via Cozzi
55, 20125 Milano, Italy. \\
Email: \url{Simone.Secchi@unimib.it} \\
URL: \url{http://www.matapp.unimib.it/~secchi}
\end{minipage}

\end{document}